\documentclass[aos,MSNbibl,seceqn,dvips]{arximspdf}

\doi{10.1214/14-AOS1259} 
\volume{43}
\issue{1}
\pubyear{2015}
\firstpage{1}
\lastpage{29}
\docsubty{FLA}

\makeatletter
\def\mid{|}
\newcommand{\rrvert}{\vert}
\newcommand{\llvert}{\vert}
\newcommand{\eqref}[1]{(\ref{#1})}
\newtheorem{thmm}{Theorem}[section]
\newtheorem{prop}[thmm]{Proposition}
\newproclaim{rmk}[thmm]{Remark}
\newproclaim{example}[thmm]{Example}
\newproclaim{defn}{Definition}[section]

\def\E{\mathbb{E}}
\def\P{\mathbb{P}}
\def\R{\mathbb{R}}
\def\N{\mathbb{N}}
\def\Z{\mathbb{Z}}

\def\cC{\mathcal{C}}

\def\sT{\mathsf{\Theta}}
\def\sX{\mathsf{X}}
\def\sY{\mathsf{Y}}

\makeatother

\begin{document}
\begin{frontmatter}

\title{Consistency of maximum likelihood estimation for some dynamical
systems}
\runtitle{Consistency of MLE for some dynamical systems}

\begin{aug}
\author[A]{\fnms{Kevin}~\snm{McGoff}\corref{}\ead[label=e1]{mcgoff@math.duke.edu}\thanksref{T1}},
\author[B]{\fnms{Sayan}~\snm{Mukherjee}\ead[label=e2]{sayan@stat.duke.edu}\thanksref{T2}},
\author[C]{\fnms{Andrew}~\snm{Nobel}\ead[label=e3]{nobel@email.unc.edu}\thanksref{T3}}
\and
\author[D]{\fnms{Natesh} \snm{Pillai}\ead[label=e4]{pillai@fas.harvard.edu}\thanksref{T4}}
\runauthor{M{\normalfont{c}}Goff, Mukherjee, Nobel and Pillai}

\thankstext{T1}{Supported by NSF Grant DMS-10-45153.}
\thankstext{T2}{Supported by NIH (Systems Biology): 5P50-GM081883,
AFOSR: FA9550-10-1-0436, NSF~CCF-1049290 and NSF DMS-12-09155.}
\thankstext{T3}{Supported in part by NSF Grants DMS-09-07177 and
DMS-13-10002.}
\thankstext{T4}{Supported by NSF Grant DMS-11-07070.}

\affiliation{Duke University, Duke University, University of North
Carolina and Harvard University}

\address[A]{K. McGoff\\
Department of Mathematics\\
Duke University\\
Durham, North Carolina 27708\\
USA\\
\printead{e1}}
\address[B]{S. Mukherjee\\
Departments of Statistical Science,\\
\quad Computer Science, and Mathematics\\
Institute for Genome Sciences \& Policy\\
Duke University\\
Durham, North Carolina 27708\\
USA\\
\printead{e2}}
\address[C]{A. Nobel\\
Department of Statistics\\
\quad and Operations Research\\
University of North Carolina \\
Chapel Hill, North Carolina 27599-3260\\
USA\\
\printead{e3}}
\address[D]{N. Pillai\\
Department of Statistics\\
Harvard University\\
Cambridge, Massachusetts 02138\hspace*{27pt}\\
USA\\
\printead{e4}}
\end{aug}
%

\received{\smonth{6} \syear{2013}}
\revised{\smonth{7} \syear{2014}}

%
\begin{abstract}
We consider the asymptotic consistency of maximum likelihood parameter
estimation for dynamical systems observed with noise. 
Under suitable conditions on the dynamical systems and the
observations, we show that maximum likelihood parameter estimation is
consistent. Our proof involves ideas from both information theory and
dynamical systems. Furthermore, we show how some well-studied
properties of dynamical systems imply the general statistical
properties related to maximum likelihood estimation. Finally, we
exhibit classical families of dynamical systems for which maximum
likelihood estimation is consistent. Examples include shifts of finite
type with Gibbs measures and Axiom~A attractors with SRB measures.
\end{abstract}



\end{frontmatter}

\section{Introduction} \label{Introduction}

Maximum likelihood estimation is a common, well-\break studied and powerful
technique for statistical estimation. In the context of a statistical
model with an unknown parameter, the maximum likelihood estimate of the
unknown parameter is, by definition, any parameter value under which
the observed data is most likely; such parameter values are said to
maximize the likelihood function with respect to the observed data. In
classical statistical models, one typically thinks of the unknown
parameter as a real number or possibly a finite dimensional vector of
real numbers. Here we consider maximum likelihood estimation for
statistical models in which each parameter value corresponds to a
stochastic system observed with noise.

Hidden Markov models (HMMs) provide a natural setting in which to study
both stochastic systems with observational noise and maximum likelihood
estimation. In this setting, one has a parametrized family of
stochastic processes that are assumed to be Markov, and one attempts to
perform inference about the underlying parameters from noisy
observations of the process. There has been a substantial amount of
work on statistical inference for HMMs, and we do not attempt a
complete survey of that area here. In the 1960s, Baum and Petrie \cite
{BaumPetrie1966,Petrie1969} studied consistency of maximum likelihood
estimation for finite state HMMs. Since that time, several other
authors have shown that maximum likelihood estimation is consistent for
HMMs under increasingly general conditions \cite
{DoucMatias2001,DMR2004,GL2006,LeGlandMevel2000,LeGlandMevel2000_2,Leroux1992},
culminating with the work of Douc et al. \cite{DMOvH2011}, which
currently provides the most general conditions on HMMs under which
maximum likelihood estimation has been shown to be consistent.

We focus here on the consistency of maximum likelihood estimation for
parametrized families of deterministic systems observed with noise.
Inference methods for deterministic systems from noisy observations are
of interest in a variety of scientific areas; for a few examples, see
\cite
{IBK2006,IBAK2011,LawStuart2012,PS2004,PooleRaftery2000,RHCC2007,TWSIS2009,Wood2010}.

For the purpose of this article, the terms deterministic system and
dynamical system refer to a map $T \dvtx\sX\to\sX$. The set $\sX$ is
referred to as the state space, and the transformation $T$ governs the
evolution of states over one (discrete) time increment. Our main
interest here lies in families of dynamical systems observed with
noise. More precisely, we consider a state space $\sX$ and a parameter
space $\sT$, and to each $\theta$ in $\sT$, we associate a dynamical
system $T_{\theta} \dvtx\sX\to\sX$. Note that the state space $\sX$ does
not depend on $\theta$. For each $\theta$ in $\sT$, we assume that the
system is started at equilibrium from a $T_{\theta}$-invariant measure
$\mu_{\theta}$. See Section~\ref{SectionAssumptions} for precise
definitions. We are particularly interested in situations in which the
family of dynamical systems is observed via noisy measurements (or
observations). We consider a general observation model specified by a
family of probability densities $\{g_{\theta}( \cdot\mid x) \dvtx\theta
\in\sT, x \in\sX\}$, where $g_{\theta}( \cdot\mid x)$ prescribes
the distribution of an observation given that the state of the
dynamical system is $x$ and the state of nature is~$\theta$. Under some
additional conditions (see Section~\ref{GeneralTheorem}), our first
main result states that maximum likelihood estimation is a consistent
method of estimation of the parameter~$\theta$.

We have chosen to state the conditions of our main consistency result
in terms of statistical properties of the family of dynamical systems
and the observations. However, these particular statistical properties
have not been directly studied in the dynamical systems literature. In
the interest of applying our general result to specific systems, we
also establish several connections between well-studied properties of
dynamical systems and the statistical properties relevant to maximum
likelihood estimation. Finally, we apply these results to some
examples, including shifts of finite type with Gibbs measures and Axiom~A attractors with SRB (Sinai--Ruelle--Bowen) measures. It is widely
accepted in the field of ergodic theory and dynamical systems that
these classes of systems have ``good'' statistical properties, and our
results may be viewed as a precise confirmation of this view.

\subsection{Previous work} \label{SectionPriorWork}

There has been a substantial amount of work on statistical inference
for HMMs, and a complete survey of that area is beyond the scope of
this work. The asymptotic consistency of maximum likelihood estimation
for HMMs has been studied at least since the work of Baum and Petrie
\cite{BaumPetrie1966,Petrie1969} under the assumption that both the
hidden state space $\sX$ and the observation space $\sY$ are finite
sets. Leroux extended this result to the setting where $\sY$ is a
general space and $\sX$ is a finite set \cite{Leroux1992}. Several
other authors have shown that maximum likelihood estimation is
consistent for HMMs under increasingly general conditions \mbox{\cite
{DoucMatias2001,DMR2004,GL2006,LeGlandMevel2000,LeGlandMevel2000_2}},
culminating with the work of Douc et al. \cite{DMOvH2011}, which
currently provides the most general conditions for HMMs under which
maximum likelihood estimation has been shown to be consistent.

Let us now discuss the results of Douc et al. \cite{DMOvH2011} in
greater detail. Consider parametrized families of HMMs in which both
the hidden state space $\sX$ and the observation space $\sY$ are
complete, separable metric spaces. The main result of \cite{DMOvH2011}
shows that under several conditions, maximum likelihood estimation is a
consistent method of estimation of the unknown parameter. These
conditions involve some requirements on the transition kernel of the
hidden Markov chain, as well as basic integrability conditions on the
observations. The proof of that result relies on information-theoretic
arguments, in combination with the application of some mixing
conditions that follow from the assumptions on the transition kernel.
To prove our consistency result, we take a similar
information-theoretic approach, but instead of placing explicit
restrictions on the transition kernel, we identify and study mixing
conditions suitable for dynamical systems. See Remarks \ref
{ComparisonRmk} and \ref{StochasticRemark} for further discussion of
our results in the context of HMMs.

Other directions of study regarding inference for HMMs include the
behavior of MLE for misspecified HMMs \cite{DoucMoulines2012},
asymptotic normality for parameter estimates \cite
{BRR1998,JensenPetersen1999}, the dynamics of Bayesian updating \cite
{Shalizi2009} and starting the hidden process away from equilibrium
\cite{DMOvH2011}. Extending these results to dynamical systems is of
potential interest.

The topic of statistical inference for dynamical systems has been
widely studied in a variety of fields. Early interest from the
statistical point of view is reflected in the following surveys \cite
{Berliner1992,ChatterjeeYilmaz1992,Isham1993,Jensen1993}. For a recent
review of this area with many references, see \cite{MMP2013}. There has
been significant methodological work in the area of statistical
inference for dynamical systems (for a few recent examples, see \cite
{IBAK2011,IBK2006,PS2004,TWSIS2009,Wood2010}), but in this section we
attempt to describe some of the more theoretical work in this area. The
relevant theoretical work to date falls (very) roughly into three classes:
\begin{itemize}
\item state estimation (also known as denoising or filtering) for
dynamical systems with observational noise;
\item prediction for dynamical systems with observational noise;
\item system reconstruction from dynamical systems without noise.
\end{itemize}
Let us now mention some representative works from these lines of research.

In the setting of dynamical systems with observational noise, Lalley
introduced several ideas regarding state estimation in \cite
{Lalley1999}. These ideas were subsequently generalized and developed
in \cite{Lalley2001,LalleyNobel2006}. Key results from this line of
study include both positive and negative results on the consistency of
denoising a dynamical system under additive observational noise. In
short, the magnitude of the support of the noise seems to determine
whether consistent denoising is possible. In related work, Judd \cite
{Judd2007} demonstrated that MLE can fail (in a particular sense) in
state estimation when noise is large. It is perhaps interesting to note
that there are examples of Axiom~A systems with Gaussian observational
noise for which state estimation cannot be consistent (by results of
\cite{Lalley2001,LalleyNobel2006}) and yet MLE provides consistent
parameter estimation (by Theorem~\ref{AxiomAThm}).

Steinwart and Anghel considered the problem of consistency in
prediction accuracy for dynamical systems with observational noise
\cite{SteinwartAnghel2009}. They were able to show that support vector
machines are consistent in terms of prediction accuracy under some
conditions on the decay of correlations of the dynamical system.

The work of Adams and Nobel uses ideas from regression to study
reconstruction of measure-preserving dynamical systems \cite
{AdamsNobel2001,Nobel2001,AdamsNobel2001_2} without noise. These
results show that certain types of inference are possible under fairly
mild ergodicity assumptions. A sample result from this line of work is
that a measure-preserving transformation may be consistently
reconstructed from a typical trajectory observed without noise,
assuming that the transformation preserves a measure that is absolutely
continuous (with Radon--Nikodym derivative bounded away from $0$ and
infinity) with respect to a known reference measure.

\subsection{Organization} \label{SectionOrganization}

In Section~\ref{SectionAssumptions}, we give some necessary background
on dynamical systems observed with noise. Section~\ref{GeneralTheorem}
contains a statement and discussion of our main result (Theorem~\ref
{mainThm}), which asserts that under some general statistical
conditions, maximum likelihood parameter estimation is consistent for
\mbox{families} of dynamical systems observed with noise. The purpose of
Section~\ref{SectionVerification} is to establish connections between
well-studied properties of dynamical systems and the (statistical)
conditions appearing in Theorem~\ref{mainThm}. Section~\ref
{SectionExamples} gives several examples of widely studied families of
dynamical systems to which we apply Theorem~\ref{mainThm} and therefore
establish consistency of maximum likelihood estimation. The proofs of
our main results appear in Section~\ref{SectionProofs}, and we conclude
with some final remarks in Section~\ref{SectionConclusion}.

\section{Setting and notation} \label{SectionAssumptions}

Recall that our primary objects of study are parametrized families of
dynamical systems. In this section we introduce these objects in some
detail. First let us recall some terminology regarding dynamical
systems and ergodic theory. We use $\sX$ to denote a state space, which
we assume to be a complete separable metric space endowed with its
Borel $\sigma$-algebra $\mathcal{X}$. Then a measurable dynamical
system on $\sX$ is defined by a measurable map $T \dvtx\sX\to\sX$, which
governs the evolution of states over one (discrete) time increment. For
a probability measure $\mu$ on the measurable space $(\sX,\mathcal
{X})$, we say that $T$ preserves $\mu$ (or $\mu$ is $T$-invariant) if
$\mu(T^{-1}E) = \mu(E)$ for each set $E$ in $\mathcal{X}$. We refer to
the quadruple $(\sX,\mathcal{X},T,\mu)$ as a measure-preserving system.
To generate a trajectory $(X_k)$ from such a measure-preserving system,
one chooses $X_0$ according to $\mu$ and sets $X_k = T^k(X_0)$ for $k
\geq0$. Note that $(X_k)$ is then a stationary $\sX$-valued stochastic
process. Finally, the measure-preserving system $(\sX,\mathcal
{X},T,\mu
)$ is said to be ergodic if $T^{-1}E = E$ implies $\mu(E) \in\{0,1\}$.
See the books \cite{Petersen,Walters} for an introduction to
measure-preserving systems and ergodic theory.

Let us now introduce the setting of parametrized families of dynamical
systems. We denote the parameter space by $\sT$, which is assumed to be
a compact metric space endowed with its Borel $\sigma$-algebra. Fix a
state space $\sX$ and its Borel $\sigma$-algebra $\mathcal{X}$ as
above. To each parameter $\theta$ in $\sT$, we associate a measurable
transformation $T_{\theta} \dvtx\sX\to\sX$, which prescribes the
dynamics corresponding to the parameter $\theta$. Finally, we need to
specify some initial conditions. In this article, we consider the case
that the system is started from equilibrium. More precisely, we
associate to each $\theta$ in $\sT$ a $T_{\theta}$-invariant Borel
probability measure $\mu_{\theta}$ on $(\sX,\mathcal{X})$. Thus, to
each $\theta$ in $\sT$, we associate a measure-preserving system
$(\sX
,\mathcal{X},T_{\theta},\mu_{\theta})$, and we refer to the collection
$(\sX,\mathcal{X},T_{\theta},\mu_{\theta})_{\theta\in\sT}$ as a
parametrized family of dynamical systems. For ease of notation, we will
refer to $(T_{\theta},\mu_{\theta})_{\theta\in\sT}$ as a family of
dynamical systems on $(\sX,\mathcal{X})$, instead of referring to the
family of quadruples $(\sX,\mathcal{X},T_{\theta},\mu_{\theta
})_{\theta
\in\sT}$.

We would like to study the situation that such a family of dynamical
systems is observed via noisy measurements. Here we describe the
specifics of our observation model. We suppose that we have a complete,
separable metric space $\sY$, endowed with its $\sigma$-algebra
$\mathcal{Y}$, which serves as our observation space. We also assume
that we have a family of Borel probability densities $\{g_{\theta
}(\cdot\mid x) \dvtx\theta\in\sT, x \in\sX\}$ with respect to a
fixed reference measure $\nu$ on $\sY$. The density $g_{\theta}(
\cdot
\mid x)$ prescribes the distribution of our observation given that the
state of the dynamical system is $x$ and the state of nature is $\theta
$. Finally, we assume that the noise involved in successive
observations is conditionally independent given $\theta$ and the
underlying trajectory of the dynamical system. Thus our full model
consists of a parametrized family of dynamical systems $(T_{\theta
},\mu
_{\theta})_{\theta\in\sT}$ on a measurable space $(\sX,\mathcal{X})$
with corresponding observation densities $\{g_{\theta}(\cdot\mid x)
\dvtx\theta\in\sT, x \in\sX\}$.

In general, we would like to estimate the parameter $\theta$ from our
observations. Maximum likelihood estimation provides a basic method for
performing such estimation. Our first main result states that maximum
likelihood estimation is a consistent estimator of $\theta$ under some
general conditions on the family of systems and the noise. In order to
state these results precisely, we now introduce the likelihood for our
model. For the sake of notation, it will be convenient to denote finite
sequences $(x_i,\ldots,x_j)$ with the notation $x_i^j$.

As we have assumed that our observations are conditionally independent
given $\theta$ and a trajectory $(X_k)$, we have that for $\theta\in
\sT$ and $y_0^n \in\sY^{n+1}$, the likelihood of observing $y_0^n$
given $\theta$ and $(X_k)$ is
\[
p_{\theta}\bigl(y_0^n \mid X_0^n
\bigr) = \prod_{j=0}^{n}
g_{\theta}(y_j \mid X_j).
\]
Since $X_k = T_{\theta}^k(X_0)$ given $\theta$ and $X_0$, the
conditional likelihood of $y_0^n$ given $\theta$ and $X_0 = x$ is
\[
p_{\theta}\bigl(y_0^n \mid x\bigr) = \prod
_{j=0}^{n} g_{\theta}
\bigl(y_j \mid T_{\theta}^j(x)\bigr).
\]
Since our model also assumes that $X_0$ is distributed according to
$\mu
_{\theta}$, we have that for $\theta\in\sT$ and $y_0^n \in\sY
^{n+1}$, the marginal likelihood of observing $y_0^n$ given $\theta$ is
%
%
\begin{equation}
\label{likelihood} p_{\theta}\bigl(y_0^n\bigr) = \int
p_{\theta}\bigl(y_0^n \mid x\bigr) \,d
\mu_{\theta}(x).
\end{equation}
We denote by $\nu^n$ the product measure on $\sY^{n+1}$ with marginals
equal to $\nu$. Let $\P_{\theta}$ be the probability measure on $\sX
\times\sY^{\N}$ such that for Borel sets $A \subset\sX$ and $B
\subset\sY^{n+1}$, it holds that
\[
\P_{\theta}(A \times B) = \int\int\mathbf{1}_A(x) \mathbf
{1}_B\bigl(y_0^n\bigr) p_{\theta}
\bigl(y_0^n \mid x\bigr) \,d\nu^{n}
\bigl(y_0^n\bigr) \,d\mu_{\theta}(x),
\]
which is well defined by Kolmogorov's consistency theorem. Let $\E
_{\theta}$ denote expectation with respect to $\P_{\theta}$, and let
$\P
_{\theta}^Y$ be the marginal of $\P_{\theta}$ on $\sY^{\N}$.

Before we define consistency, let us first consider the issue of
identifiability. Our notion of identifiability is captured by the
following equivalence relation.
%

\begin{defn}
Define an equivalence relation on $\sT$ as follows: let $\theta\sim
\theta'$ if $\P_{\theta}^Y = \P_{\theta'}^Y$. Denote by $[\theta
]$ the
equivalence class of $\theta$ with respect to this equivalence relation.
\end{defn}

In a strong theoretical sense, if $\theta'$ is in $[\theta]$, then the
systems corresponding to the parameter values $\theta'$ and $\theta$
cannot be distinguished from each other based on observations of the system.

Now we fix a distinguished element $\theta_0$ in $\sT$. Here and in the
rest of the paper, we assume that $\theta_0$ is the ``true'' parameter;
that is, the data are generated from the measure $\mathbb{P}^Y_{\theta_0}$.
Hence, one may think of $[\theta_0]$ as the set of parameters that
cannot be distinguished from the true parameter.

%
\begin{defn}
An approximate maximum likelihood estimator (MLE) is a sequence of
measurable functions $\hat{\theta}_n \dvtx(\sY)^{n+1} \to\sT$ such that
%
%
\begin{equation}
\label{MLEdef} \frac{1}{n} \log p_{\hat{\theta}_n(Y_0^{n})}\bigl(Y_0^{n}
\bigr) \geq\sup_{\theta
} \frac{1}{n} \log p_{\theta}
\bigl(Y_0^{n}\bigr) - o_{\mathrm{a.s.}}(1),
\end{equation}
where $o_{\mathrm{a.s.}}(1)$ denotes a process that tends to zero $\P_{\theta
_0}$-a.s. as $n$ tends to infinity.
\end{defn}

%
\begin{rmk} Several notions in this article, including the definition
of approximate MLE above, involve taking suprema over $\theta$ in $\sT
$. In many situations of interest to us, $\sX$ and $\sT$ are compact,
and all relevant functions are continuous in these arguments. In such
cases, we have sufficient regularity to guarantee that suprema over
$\theta$ in $\sT$ are measurable. However, in the general situation,
such suprema are not guaranteed to be measurable, and one must take
some care. As all our measurable spaces are Polish (complete, separable
metric spaces); such functions are always universally measurable \cite
{BertsekasShreve1978}, Proposition~7.47. Similarly, a~Borel-measurable
(approximate) maximum likelihood estimator need not exist, but the
Polish assumption ensures the existence of universally measurable
maximum likelihood estimators \cite{BertsekasShreve1978},
Proposition~7.50. Thus all probabilities and expectations may
be unambiguously extended to such quantities.
\end{rmk}

%
\begin{rmk}
In this work, we do not consider specific schemes for constructing an
approximate MLE. Based on the existing results regarding denoising and
system reconstruction (e.g., \cite
{AdamsNobel2001,Lalley1999,Lalley2001,LalleyNobel2006,Nobel2001,AdamsNobel2001_2},
which are briefly discussed in Section~\ref{SectionPriorWork}),
explicit construction of an approximate MLE may be possible under
suitable conditions. Although the description and study of such
constructive methods could be interesting, it is outside of the scope
of this work.
\end{rmk}

%
\begin{rmk}
In principle, one could consider inference based on the conditional
likelihood $p_{\theta}( \cdot\mid x_0)$ in place of the marginal
likelihood $p_{\theta}( \cdot)$. However, we do not pursue this
direction in this work. For nonlinear dynamical systems, even the
conditional likelihood $p_{\theta}( \cdot\mid x_0)$ may depend very
sensitively on $x_0$; see~\cite{Berliner1992}, for example. Thus
optimizing over $x_0$ is essentially no more ``tractable'' than
marginalizing the likelihood via an invariant measure.
\end{rmk}

%
\begin{rmk} \label{ComparisonRmk}
The framework of this paper may be translated into the language of
Markov chains as follows. For each $\theta\in\sT$, we define a
(degenerate) Markov transition kernel $Q_{\theta}$ as follows:
\[
Q_{\theta}(x,y) = \delta_{T_{\theta}(x)}(y).
\]
In other words, for each $\theta\in\sT$, $x \in\sX$, and Borel set
$A \subset\sX$, the probability that $X_1 \in A$ conditioned on $X_0 =
x$ is
\[
Q_{\theta}(x,A) = \delta_{T_{\theta}(x)}(A),
\]
where $\delta_x$ is defined to be a point mass at $x$.

In all previous work on consistency of maximum likelihood estimation
for HMMs (including \cite
{DoucMatias2001,DMOvH2011,DMR2004,GL2006,LeGlandMevel2000,LeGlandMevel2000_2}),
there have been significant assumptions placed on the Markovian
structure of the hidden chain. For example, the central hypothesis
appearing in \cite{DMOvH2011} requires that there is a $\sigma$-finite
measure $\lambda$ on $\sX$ such that for some $L \geq0$, the $L$-step
transition kernel $Q_{\theta}^L(x,\cdot)$ is absolutely continuous with
respect to $\lambda$ with bounded Radon--Nikodym derivative. If $\sX$
is uncountable, then the degeneracy of $Q_{\theta}$, which arises
directly from the fact that we are considering deterministic systems,
makes the existence of such a dominating measure impossible. In short,
it is precisely the determinism in our hidden processes that prevents
previous theorems for HMMs from applying to dynamical systems.

Nonetheless, there is a special case of systems that we consider in
Section~\ref{SectionGibbs} that overlaps with the systems considered in
the HMM literature. If $\sX$ is a shift of finite type, $T_{\theta}$ is
the shift map $\sigma\dvtx \sX\to\sX$ for all $\theta$, $\mu_{\theta}$
is a ($1$-step) Markov measure for all $\theta$, and $g_{\theta}(
\cdot
\mid x)$ depends only $\theta$ and the zero coordinate $x_0$, then both
the present work and the results in \cite{DMOvH2011} apply to this
setting and guarantee consistency of any approximate MLE under
additional assumptions on the noise.
\end{rmk}

\section{Consistency of MLE} \label{GeneralTheorem}

In this section, we show that under suitable conditions, any
approximate MLE is consistent for families of dynamical systems
observed with noise. To make this statement precise, we make the
following definition of consistency.

%
\begin{defn}
An approximate MLE $(\hat{\theta}_n)_n$ is consistent at $\theta_0$ if
$\hat{\theta}_n(Y_0^n)$ converges to $[\theta_0]$, $\P_{\theta_0}$-a.s.
as $n$ tends to infinity.
\end{defn}

For the sake of notation, define the function $\gamma\dvtx \sT\times\sY
\to\R_+$, where
\[
\gamma_{\theta}(y) = \sup_{x \in\sX} g_{\theta}(y \mid
x).
\]
Also, for $x >0$, let $\log^+ x = \max(0,\log(x))$.

Consider the following conditions on a family of dynamical systems
observed with noise:
\begin{longlist}[(S1)]
\item[(S1)]\textit{Ergodicity}. 

The system $(T_{\theta_0},\mu_{\theta_0})$ on $(\sX
,\mathcal
{X})$ is ergodic.

\item[(S2)]\textit{Logarithmic integrability at} $\theta_0$. 

It holds that
\[
\E_{\theta_0} \bigl[ \log^+ \gamma_{\theta_0}(Y_0) \bigr] <
\infty
\]
and
\[
\E_{\theta_0} \biggl[ \biggl\llvert\log\int g_{\theta_0}(Y_0
\mid x) \,d\mu_{\theta_0}(x) \biggr\rrvert\biggr] < \infty.
\]

\item[(S3)]\textit{Logarithmic integrability away from $\theta_0$.} 

For each $\theta' \notin[\theta_0]$, there exists a
neighborhood $U$ of $\theta'$ such that
\[
\E_{\theta_0} \Bigl[ \sup_{\theta\in U} \log^+
\gamma_{\theta}(Y_0) \Bigr] < \infty.
\]

\item[(S4)]\textit{Upper semi-continuity of the likelihood}. 

For each $\theta' \notin[\theta_0]$ and $n \geq0$, the
function $\theta\mapsto p_{\theta}(Y_0^n)$ is upper semi-continuous at
$\theta'$, $\P_{\theta_0}$-a.s.

\item[(S5)]\textit{Mixing condition}. 

There exists $\ell\geq0$ such that for each $m \geq0$,
there exists a measurable function $C_{m} \dvtx\sT\times\sY^{m+1} \to
\R
_+$ such that if $t \geq1$ and $w_0, \ldots, w_t \in\sY^{m+1}$, then
\[
\int\prod_{j=0}^t p_{\theta}
\bigl(w_j \mid T_{\theta}^{j(m+\ell)}x\bigr) \,d\mu
_{\theta}(x) \leq\prod_{j=0}^t
C_{m}(\theta,w_j) \prod_{j=0}^t
p_{\theta}(w_j).
\]
Furthermore, for each $\theta' \notin[\theta_0]$, there exists a
neighborhood $U$ of $\theta'$ such that
\[
\sup_m \E_{\theta_0} \Bigl[ \sup_{\theta\in U}
\log C_{m}\bigl(\theta,Y_0^{m}\bigr) \Bigr] <
\infty.
\]

\item[(S6)]\textit{Exponential identifiability}. 

For each $\theta\notin[\theta_0]$, there exists a sequence
of measurable sets $A_n \subset\sY^{n+1}$ such that
\[
\liminf_n \mathbb{P}^Y_{\theta_0}(A_n)
>0 \quad\mbox{and}\quad \limsup_n \frac
{1}{n} \log
\mathbb{P}^Y_{\theta}(A_n) < 0.
\]
\end{longlist}

The following theorem is our main general result.

%
\begin{thmm} \label{mainThm}
Suppose that $(T_{\theta},\mu_{\theta})_{\theta\in\sT}$ is a
parametrized family of dynamical systems on $(\sX,\mathcal{X})$ with
corresponding observation densities $(g_{\theta})_{\theta\in\sT}$. If
conditions \textup{(S1)--(S6)} hold, then any
approximate MLE is consistent at $\theta_0$.
\end{thmm}

The proof of Theorem~\ref{mainThm} is given in Section~\ref
{SectionProofs}. In the following remark, we discuss conditions (S1)--(S6).

%
\begin{rmk} \label{ConditionsRmk}
Conditions (S1)--(S3) involve basic
irreducibility and integrability conditions, and similar conditions
have appeared in previous work on consistency of maximum likelihood
estimation for HMMs; see, for example, \mbox{\cite{DMOvH2011,Leroux1992}}.
Taken together, conditions (S1) and (S2) ensure the almost sure existence and finiteness of
the entropy rate for the process $(Y_n)$,
\[
h(\theta_0) = \lim_n \frac{1}{n} \log
p_{\theta_0}\bigl(Y_0^n\bigr).
\]
Condition (S3) serves as a basic integrability
condition in the proof of Theorem~\ref{mainThm}, in which one must
essentially show that for $\theta\notin[\theta_0]$,
\[
\limsup_n \frac{1}{n} \log p_{\theta}
\bigl(Y_0^n\bigr) < h(\theta_0).
\]
Conditions (S4)--(S6) are more interesting
from the point of view of dynamical systems, and we discuss them in
greater detail below.

The upper semi-continuity of the likelihood (S4) is closely
related to the continuity of the map $\theta\mapsto\mu_{\theta}$. In
general, the continuous dependence of $\mu_{\theta}$ on $\theta$ places
nontrivial restrictions on a family of dynamical systems. This
property (continuity of $\theta\mapsto\mu_{\theta}$) is often called
``statistical stability'' in the dynamical systems and ergodic theory
literature, and it has been studied for some families of systems; for example,
see \cite{ACF2010,FreitasTodd2009,Ruelle1997,Vasquez2007} and
references therein. In Section~\ref{StatStability}, we show how
statistical stability of the family of dynamical systems may be used to
establish the upper semi-continuity of the likelihood (S4).

The mixing condition (S5) involves control of the
correlations of the observation densities along trajectories of the
underlying dynamical system. Although the general topic of decay of
correlations has been widely studied in dynamical systems (see \cite
{Baladi2001} for an overview), condition (S5) is not
implied by the particular decay of correlations properties that are
typically studied for dynamical systems. Nonetheless, we show in
Section~\ref{SectionSubadd} how some well-studied mixing properties of
dynamical systems imply the mixing condition (S5).

Finally, condition (S6) involves the exponential
identifiability of the true parameter $\theta_0$. We show in
Section~\ref{SectionIdentifiability} how large deviations for a family of
dynamical systems may be used to establish exponential identifiability
(S6). Large deviations estimates for dynamical
systems have been studied in \cite{RBY2008,Young1990}, and our main
goal in Section~\ref{SectionIdentifiability} is to connect such results
to exponential identifiability (S6).
\end{rmk}

%
\begin{rmk} \label{StochasticRemark}
Suppose one has a family of bi-variate stochastic processes $\{
(X_k^{\theta},Y_k^{\theta}) \dvtx\theta\in\sT\}$, where $(X_k^{\theta})$
is interpreted as a hidden process and $(Y_k^{\theta})$ as an
observation process. If the observations have conditional densities
with respect to a common measure given $(X_k^{\theta})$ and $\theta$,
then it makes sense to ask whether maximum likelihood estimation is a
consistent method of inference for the parameter $\theta$.

It is well known that the setting of stationary stochastic processes
may be translated into the deterministic setting of dynamical systems,
which may be carried out as follows. Let $\{(X_k^{\theta}) \dvtx\theta
\in
\sT\}$ be a family of stationary stochastic processes on a measurable
space $(\sX,\mathcal{X})$. Consider the product space $\hat{\sX} =
\sX
^{\otimes\Z}$ with corresponding $\sigma$-algebra $\hat{\mathcal{X}}$.
Each process $(X_k^{\theta})$ corresponds to a probability measure
$\mu
_{\theta}$ on $(\hat{\sX},\hat{\mathcal{X}})$ with the property that
$\mu_{\theta}$ is invariant under the left-shift map $T\dvtx\hat{\sX}
\to
\hat{\sX}$ given by $\mathbf{x} = (x_i)_i \mapsto T(\mathbf{x}) =
(x_{i+1})_i$. With this translation, Theorem~\ref{mainThm} shows that
maximum likelihood estimation is consistent for families of hidden
stochastic processes $(X_k^{\theta})$ observed with noise, whenever the
corresponding family of dynamical systems $(T,\mu_{\theta})$ on
$(\hat
{\sX},\hat{\mathcal{X}})$ with observation densities satisfy conditions
(S1)--(S6).

With the above translation, Theorem~\ref{mainThm} applies to some
families of processes allowing infinite-range dependence in both the
hidden process $(X_k^{\theta})$ and the observation process
$(Y_k^{\theta})$. From this point of view, Theorem~\ref{mainThm}
highlights the fact that maximum likelihood estimation is consistent
for dependent processes observed with noise as long as they satisfy
some general conditions: ergodicity, logarithmic integrability of
observations, continuous dependence on the parameters and some mixing
of the observation process. It is interesting to note that the existing
work on consistency of maximum likelihood estimation for HMMs \cite
{CMR2005,DoucMatias2001,DMOvH2011,DMR2004,GL2006,LeGlandMevel2000,LeGlandMevel2000_2,Leroux1992}
makes assumptions of precisely this sort in the specific context of
Markov chains.
\end{rmk}

\section{Statistical properties of dynamical systems} \label
{SectionVerification}

In our main consistency result (Theorem~\ref{mainThm}), we establish
the consistency of any approximate MLE under conditions (S1)--(S6). We have chosen to formulate
our result in these terms because they reflect general statistical
properties of dynamical systems observed with noise that are relevant
to parameter inference. However, these conditions have not been
explicitly studied in the dynamical systems literature, despite the
fact that much effort has been devoted to understanding certain
statistical aspects of dynamical systems. In this section, we make
connections between the general statistical conditions appearing in
Theorem~\ref{mainThm} and some well-studied properties of dynamical
systems. Section~\ref{StatStability} shows how the notion of
statistical stability may be used to verify the upper semi-continuity
of the likelihood (S4). Section~\ref{SectionSubadd} connects
well-known mixing properties of some measure-preserving dynamical
systems to the mixing property (S5). In Section~\ref
{SectionIdentifiability}, we show how large deviations for dynamical
systems may be used to deduce the exponential identifiability condition~(S6). Proofs of statements in this section, as well
as additional discussion, appear in Supplementary Appendix~A \cite{supplement}.

\subsection{Statistical stability and continuity of \texorpdfstring{$p_{\theta}$}{ptheta}}
\label
{StatStability}

As discussed in Remark~\ref{ConditionsRmk}, the upper semi-continuity
condition (S4) places nontrivial restrictions on the family of
dynamical systems under consideration. In this section, we establish
sufficient conditions for (S4) to hold. The continuous
dependence of $\mu_{\theta}$ on $\theta$ is a property called
statistical stability in the dynamical systems literature \cite
{ACF2010,FreitasTodd2009,Ruelle1997,Vasquez2007}. Let us state this
property precisely. Let $M(\sX)$ denote the space of Borel probability
measures on $\sX$. Endow $M(\sX)$ with the topology of weak
convergence: $\mu_n$ converges to $\mu$ if $\int f \,d\mu_n$ converges to
$\int f \,d\mu$ as $n$ tends to infinity, for each continuous, bounded
function $f\dvtx\sX\to\R$. The family of dynamical systems $(T_{\theta
},\mu_{\theta})_{\theta\in\sX}$ on $(\sX,\mathcal{X})$ is said to
have statistical stability if the map $\theta\mapsto\mu_{\theta}$ is
continuous with respect to the weak topology on $M(\sX)$.

The following proposition shows that under some continuity and
compactness assumptions, statistical stability of the family of
dynamical systems implies upper semi-continuity of the likelihood
(S4).
%

\begin{prop} \label{USClemma2}
Suppose that $\sX$ and $\sT$ are compact, and the maps $T \dvtx\sT
\times
\sX\to\sX$ and $g \dvtx\sT\times\sX\times\sY\to\R_+$ are
continuous. If the family $(T_{\theta},\mu_{\theta})_{\theta\in\sT}$
has statistical stability, then upper semi-continuity of the likelihood
\textup{(S4)} holds.
\end{prop}

The proof of Proposition~\ref{USClemma2} appears in Supplementary
Appendix A.1 \cite{supplement}.

\subsection{Mixing} \label{SectionSubadd}

In this section, we focus on mixing condition (S5).
Recall that (S5) involves a nontrivial restriction
on the correlations of the observation densities $g_{\theta}$ along
trajectories of the underlying dynamical system. 
Although mixing conditions have been widely studied in the dynamics
literature, the particular type of condition appearing in (S5) appears not to have been investigated. Nonetheless, we
show that a well-studied mixing property for dynamical systems implies
the statistical mixing property (S5).

In order to study mixing for dynamical systems, one typically places
restrictions on the type of events or observations that one considers
(by considering certain functionals of the process). For example, in
some situations a substantial amount work has been devoted to finding
particular partitions of state space with respect to which the system
possess good mixing properties; an example of such partitions are the
well-known Markov partitions \cite{Bowen1970}. If a system has good
mixing properties with respect to a particular partition, and if that
partition possesses certain (topological) regularity properties, then
it is often possible to show that the system also has good mixing
properties for related function classes, such as Lipschitz or H\"
{o}lder continuous observables. For variations of this approach to
mixing in dynamical systems, see the vast literature on decay of
correlations; for an introduction, see the survey \cite{Baladi2001}.

In this section, we follow the above approach to study mixing condition
(S5) for dynamical systems observed with noise.
First, we define a mixing property for families of dynamical systems
with respect to a partition (M1). Second, we
define a regularity property for partitions (M2).
Third, we define a topological regularity property for a family of
observation densities (M3). Finally, in the main
result of this section (Proposition~\ref{SubaddPropMain}), we show how
these three properties together imply the mixing condition (S5).

Here and in the rest of this section, we consider only invertible
transformations. It is certainly possible to modify the definitions
slightly to handle the noninvertible case, but we omit such modifications.

We will have need to consider finite partitions of $\sX$. The join of
two partitions $\cC_0$ and $\cC_1$ is defined to be the common
refinement of $\cC_0$ and $\cC_1$, and it is denoted $\cC_0 \vee
\cC
_1$. Note that for any measurable transformation $T \dvtx\sX\to\sX$, if
$\cC$ is a partition, then so is $T^{-1}\cC= \{T^{-1}A \dvtx A \in\cC\}$.
For a fixed partition $\cC$ and $i \leq j$, let $\cC_i^j = \bigvee
_{k=i}^j T_{\theta}^{-k}\cC$. Notice that $\cC_i^j$ depends on
$\theta$
through $T_{\theta}$, although we suppress this dependence in our
notation. Now consider the following alternative conditions, which may
be used in place of condition (S5):
\begin{longlist}[(M1)]

\item[(M1)]\textit{Mixing condition with respect to the partition} $\cC$.

There exists $L \dvtx\sT\to\R_+$ and $\ell\geq0$ such that
for all $\theta\in\sT$, $m,n \geq0$, $A \in\cC_0^m$ and $B \in
\cC
_{0}^{n}$, it holds that
\[
\mu_{\theta} \bigl(A \cap T_{\theta}^{-(m+\ell)}B \bigr) \leq
L_{\theta
} \mu_{\theta}(A) \mu_{\theta}(B).
\]
Furthermore, for each $\theta' \notin[\theta_0]$ there exists a
neighborhood $U$ of $\theta'$ such that
\[
\sup_{\theta\in U} L_{\theta} < \infty.
\]

\item[(M2)]\textit{Regularity of the partition} $\cC$. 
There exists $\beta\in(0,1)$ such that for all $\theta\in\sT$ and
$m,n \geq0$, if $A \in\cC_{-m}^n$ and $x,z \in A$, then
\[
d(x,z) \leq\beta^{\min(m,n)}.
\]

\item[(M3)]\textit{Regularity of observations}. 
There exists a function $K \dvtx\sT\times\sY\to\R_+$ such that for $y
\in\sY$ and $x,z \in\sX$,
\[
g_{\theta}(y \mid x) \leq g_{\theta}(y \mid z) \exp\bigl( K(\theta
,y) d(x,z) \bigr).
\]
Furthermore, for each $\theta' \notin[\theta_0]$, there exists a
neighborhood $U$ of $\theta'$ such that
\[
\E_{\theta_0} \Bigl[ \sup_{\theta\in U} K(\theta,Y_0)
\Bigr] < \infty.
\]
\end{longlist}

Let us now state the main proposition of this section, whose proof is
deferred to Supplementary Appendix~A.2 \cite
{supplement}.
%

\begin{prop} \label{SubaddPropMain}
Suppose $(T_{\theta},\mu_{\theta})_{\theta\in\sT}$ is a family of
dynamical systems on $(\sX,\mathcal{X})$ with corresponding observation
densities $(g_{\theta})_{\theta\in\sT}$. If there exists a partition
$\cC$ of $\sX$ such that conditions \textup{(M1)} and
\textup{(M2)} are satisfied, and if the observation regularity
condition \textup{(M3)} is satisfied, then mixing property
\textup{(S5)} holds.
\end{prop}

\subsection{Exponential identifiability} \label{SectionIdentifiability}

In this section, we study exponential identifiability condition (S6). We show how large deviations for dynamical systems
may be used in combination with some regularity of the observation
densities to establish exponential identifiability (S6).

Let $\sX_1$ and $\sX_2$ be metric spaces with metrics $d_1$ and $d_2$,
respectively. Recall that a function $f\dvtx\sX_1 \to\sX_2$ is said to be
H\"{o}lder continuous if there exist $\alpha>0$ and $C>0$ such that
for each $x,z$ in $\sX_1$, it holds that
\[
d_2\bigl(f(x),f(z)\bigr) \leq C d_1(x,z)^{\alpha}.
\]
If $(T,\mu)$ is a dynamical system on $(\sX,\mathcal{X})$ such that $T
\dvtx\sX\to\sX$ is H\"{o}lder continuous, then we refer to $(T,\mu)$ as
a H\"{o}lder continuous dynamical system. For many dynamical systems,
the class of H\"{o}lder continuous functions $f \dvtx\sX\to\R$ provides
a natural class of observables whose statistical properties are fairly
well understood and satisfy some large deviations estimates \cite
{RBY2008,Young1990}.

Consider the following conditions, which we later show are sufficient
to guarantee exponential identifiability (S6):
\begin{longlist}[(L1)]
\item[(L1)]\textit{Large deviations.} 
For each $\theta\notin[\theta_0]$, for each H\"{o}lder continuous
function $f \dvtx\sX\to\R$, and for each $\delta>0$, it holds that
\[
\limsup_n \frac{1}{n} \log\mu_{\theta} \Biggl(
\Biggl| \frac
{1}{n} \sum_{k=0}^{n-1} f\bigl(
T_{\theta}^k(x)\bigr) - \int f \,d\mu_{\theta} \Biggr| > \delta
\Biggr) < 0.
\]

\item[(L2)]\textit{Regularity of observations.} \label{HolderContinuityCondition}
There exists $\alpha>0$ and $K \dvtx\sT\times\sY\to\R_+$ such that
for each $x$ and $z$ in $\sX$, it holds that
\[
g_{\theta}(y \mid x) \leq g_{\theta}(y \mid z) \exp\bigl( K(
\theta,y) d(x,z)^{\alpha} \bigr).
\]
Furthermore, for $\theta\in\sT$ and $C>0$, it holds that
\[
\sup_x \int\exp\bigl( C K(\theta,y) \bigr)
g_{\theta}(y \mid x) \,d\nu(y) < \infty.
\]
\end{longlist}

The following proposition relates large deviations for dynamical
systems to the exponential identifiability condition (S6).

%
\begin{prop} \label{IdentifiabilityFromLargeDeviations}
Suppose that $(T_{\theta},\mu_{\theta})_{\theta\in\sT}$ is a family
of H\"{o}lder continuous dynamical systems on the $(\sX,\mathcal{X})$
with corresponding observation densities $(g_{\theta})_{\theta\in\sT
}$. Further suppose that the large deviations property \textup{(L1)} and the
observation regularity property
\textup{(L2)} are satisfied. Then the exponential
identifiability condition \textup{(S6)} holds.
\end{prop}

The proof of Proposition~\ref{IdentifiabilityFromLargeDeviations}
appears in Supplementary Appendix~A.3~\cite{supplement}.

\section{Examples} \label{SectionExamples}

In this section we present some classical families of dynamical systems
for which maximum likelihood estimation is consistent. We begin in
Section~\ref{SectionGibbs} by considering symbolic dynamical systems
called shifts of finite type. The state space for such systems consists
of (bi-)infinite sequences of symbols from a finite set, and the
transformation on the state space is always given by the ``left-shift''
map, which just shifts each point one coordinate to the left. Such
systems are considered models of ``chaotic'' dynamical systems that may
be defined by a finite amount of combinatorial information. In this
setting Gibbs measures form a natural class of invariant measures,
which have been studied due to their connections to statistical
physics. These measures play a central role in a topic called the
\textit{thermodynamic formalism}, which is well described in the books
\cite{Bowen,Ruelle}. Note that $k$th order finite state Markov chains
form a
special case of Gibbs measures. The main result of this section is
Theorem~\ref{GibbsFamiliesThm}, which states that under sufficient
regularity conditions, any approximate maximum likelihood estimator is
consistent for families of Gibbs measures on a shift of finite type.
The crucial assumptions for this theorem involve continuous dependence
of the Gibbs measures on $\theta$ and sufficiently regular dependence
of $g_{\theta}(y \mid x)$ on $x$. Additional proofs and discussion for
this section appear in the Supplementary Appendix~B
\cite{supplement}.

Having established consistency of maximum likelihood estimation for
families of Gibbs measures on a shift of finite type, we deduce in
Section~\ref{SectionFactors} that maximum likelihood estimation is
consistent for families of Axiom~A attractors observed with noise.
Axiom~A systems are well studied differentiable dynamical systems on
manifolds that, like shifts of finite type, exhibit ``chaotic''
behavior; for a thorough treatment of Axiom~A systems, see the book
\cite{Bowen}. In related statistical work, Lalley~\cite{Lalley1999}
considered the problem of denoising the trajectories of Axiom~A
systems. For these systems, there is a natural class of measures, known
as SRB (Sinai--Ruelle--Bowen) measures. See the article \cite
{Young2002} for an introduction to these measures with discussion of
their interpretation and importance. With the construction of Markov
partitions \cite{Bowen1970,Bowen}, one may view an Axiom~A attractor
with its SRB measure as a factor of a shift of finite type with a Gibbs
measure. Using this natural factor structure, we establish the
consistency of any approximate maximum likelihood estimator for Axiom~A
systems. Proofs and discussion of these topics appear in the
Supplementary Appendix~C \cite{supplement}.

\subsection{Gibbs measures} \label{SectionGibbs}

In this section, we consider the setting of symbolic dynamics, shifts
of finite type and Gibbs measures. We prove that any approximate
maximum likelihood estimator is consistent for these systems
(Theorem~\ref{GibbsFamiliesThm}) under some general assumptions on the
observations. Finally, we consider two examples of observations in
greater detail. In the first example, we consider ``discrete''
observations, corresponding to a ``noisy channel.'' In the second
example, we consider making real-valued observations with Gaussian
observational noise. For a brief introduction to shifts of finite type
and Gibbs measures that contains everything needed in this work, see
the Supplementary Appendix~B \cite{supplement}. For
a complete introduction to shifts of finite type and Gibbs measures,
see \cite{Bowen}.

Let us now consider some families of measure-preserving systems on
SFTs. Let $A$ be an alphabet, and let $M$ be a binary matrix with
dimensions $|A|\times|A|$. Let $\sX= X_M$ be the associated SFT, and
let $\mathcal{X}$ be the Borel $\sigma$-algebra on $\sX$. For
$\alpha
>0$, let $f \dvtx\sT\to C^{\alpha}(\sX)$ be a continuous map, and let
$\mu
_{\theta}$ be the Gibbs measure associated to the potential function
$f_{\theta}$. In this setting, we refer to $(\mu_{\theta})_{\theta
\in
\sT}$ as a continuously parametrized family of Gibbs measures on $(\sX
,\mathcal{X})$.
%

\begin{thmm} \label{GibbsFamiliesThm}
Suppose $\sX= X_M$ is a mixing shift of finite type and $(\mu_{\theta
})_{\theta\in\sT}$ is a continuously parametrized family of Gibbs
measures on $(\sX,\mathcal{X})$. If the family of observation densities
$(g_{\theta})_{\theta\in\sT}$ satisfies the integrability conditions
\textup{(S2)} and~\textup{(S3)} and the regularity
conditions \textup{(M3)} and \textup{(L2)}, then any approximate maximum likelihood
estimator is consistent.
\end{thmm}

The proof of Theorem~\ref{GibbsFamiliesThm} is based on an appeal to
Theorem~\ref{mainThm}. However, in order to verify the hypotheses of
Theorem~\ref{mainThm}, we combine the results of Section~\ref
{SectionVerification} with some well-known properties of Gibbs
measures. This proof appears in the Supplementary Appendix~B \cite{supplement}.

%

%
\begin{rmk} \label{OneSidedRmk}
There is an analogous\vspace*{1pt} theory of ``one-sided'' symbolic dynamics and
Gibbs measures, in which $A^{\Z}$ is replaced by $A^{\N}$ and
appropriate modifications are made in the definitions. The two-sided
case deals with invertible dynamical systems, whereas the one-sided
case handles noninvertible systems. We have stated Theorem~\ref{GibbsFamiliesThm} in the invertible setting, although it applies
as well in the noninvertible setting, with the obvious modifications.
\end{rmk}

%
\begin{example}
In this example, we consider families of dynamical systems $(T_{\theta
},\mu_{\theta})$ on $(\sX,\mathcal{X})$, where $\sX$ is a mixing shift
of finite type, $T_{\theta} = \sigma|_{\sX}$, and $\mu_{\theta}$
is a
continuous family of Gibbs measures on $\sX$ (as in Theorem~\ref
{GibbsFamiliesThm}). Here we consider the particular observation model
in which our observations of $\sX$ are passed through a discrete,
memoryless, noisy channel. Suppose that $\sY$ is a finite set, $\nu$ is
counting measure on $\sY$ and for each symbol $a$ in $A$ and parameter
$\theta$ in $\sT$, we have a probability distribution $\pi_{\theta}(
\cdot\mid a)$ on $\sY$. We consider the case that our observation
densities $g_{\theta}$ satisfy $g_{\theta}(\cdot\mid x) = \pi
_{\theta
}( \cdot\mid x_0)$. This situation is covered by Theorem~\ref
{GibbsFamiliesThm}, since the following conditions may be easily
verified: observation integrability (S2) and (S3) and observation regularity (M3)
and (L2).
\end{example}

%
\begin{example} \label{GibbsExample}
In this example, we once again consider families of dynamical systems
$(T_{\theta},\mu_{\theta})$ on $(\sX,\mathcal{X})$, such that $\sX
$ is
a mixing shift of finite type, $T_{\theta} = \sigma|_{\sX}$ and $\mu
_{\theta}$ is a continuous family of Gibbs measures on $\sX$ (as in
Theorem~\ref{GibbsFamiliesThm}). Here we consider the particular
observation model in which we make real-valued, parameter-dependent
measurements of the system, which are corrupted by Gaussian noise with
parameter-dependent variance. More precisely, let us assume that $\sY=
\R$, and there exists a Lipschitz continuous $\varphi\dvtx \sT\times\sX
\to\R$ and continuous $s \dvtx\sT\to(0,\infty)$ such that
\[
g_{\theta}(y \mid x) = \frac{1}{s(\theta)\sqrt{2\pi}} \exp\biggl( -
\frac{1}{2s(\theta)^2}
\bigl(\varphi_{\theta}(x)-y\bigr)^2 \biggr).
\]

We now proceed to verify conditions (S2), (S3), (M3) and (L2). First, by compactness and continuity,
there exist $C_1,C_2,C_3>0$ such that for $\theta$ in $\sT$, $y$ in
$\sY
$ and $x$ in $\sX$, it holds that
%
%
\begin{equation}
\label{NoiseIntEqn} C_1^{-1} \exp\bigl(-C_2y^2
\bigr) \leq g_{\theta}(y \mid x) \leq C_1 \exp\bigl( -
C_3 y^2\bigr).
\end{equation}
From (\ref{NoiseIntEqn}), one easily obtains the observation
integrability conditions (S2) and~(S3). Furthermore, there exists $C_4,C_5>0$ such that for
$x,z \in\sX$, it holds that
%
%
\begin{eqnarray}
\label{NoiseRegEqn} %
&&\frac{g_{\theta}(y \mid x)}{g_{\theta}(y \mid z)} \nonumber\\
&&\qquad = \exp\biggl(
-\frac
{1}{2 s(\theta)^2}
\bigl[ \bigl(\varphi_{\theta}(x)-y\bigr)^2 - \bigl(
\varphi_{\theta
}(z)-y\bigr)^2 \bigr] \biggr)
\nonumber
\\[-8pt]
\\[-8pt]
\nonumber
&&\qquad = \exp\biggl( -\frac{1}{2 s(\theta)^2} \bigl[ \bigl(\varphi_{\theta
}(x) -
\varphi_{\theta}(z)\bigr) \bigl(\varphi_{\theta}(x) +
\varphi_{\theta}(z)\bigr) + 2 y \bigl(\varphi_{\theta}(z) -
\varphi_{\theta}(x)\bigr)\bigr] \biggr)
\\
&&\qquad \leq\exp\bigl( \bigl(C_4 + C_5|y|\bigr)\bigl|\varphi_{\theta}(x)-
\varphi_{\theta
}(z)\bigr| \bigr).
\nonumber
\end{eqnarray}
Let $\varphi$ be Lipschitz continuous with constant $C_6$, and let
$K(\theta,y) = C_6(C_4+C_5|y|)$. With this choice of $K$ and (\ref
{NoiseRegEqn}), one may easily verify the observation regularity
conditions (M3) and (L2).\vadjust{\goodbreak}
\end{example}

%
\begin{rmk}
Similar calculations to those in Example~\ref{GibbsExample} imply that
any approximate maximum likelihood estimator is also consistent if the
observational noise is ``double-exponential'' [i.e., $g_{\theta}(y
\mid
x) \propto e^{- | y-x|}$]. Indeed, these calculations should hold for
most members of the exponential family, although we do not pursue them here.
\end{rmk}

\subsection{Axiom~A systems} \label{SectionFactors}

In this section, we show how the previous results may be applied to
some smooth (differentiable) families of dynamical systems. These
results follow easily from the results in Section~\ref{SectionGibbs},
using the work of Bowen and others (see \cite{Bowen1970,Bowen} and
references therein) in constructing Markov partitions for these
systems. With Markov partitions, Axiom~A systems may be viewed as
factors of the shifts of finite type with Gibbs measures. For a brief
introduction of Axiom~A systems that contains the details necessary for
this work, see the Supplementary Appendix~C \cite
{supplement}.

The basic fact that allows us to transfer our results from shifts of
finite type to Axiom~A systems is that consistency of maximum
likelihood estimation is preserved under taking appropriate factors.
Let us now make this statement precisely. Suppose that $(T_{\theta
},\mu
_{\theta})_{\theta\in\sT}$ is a family of dynamical systems on
$(\sX
,\mathcal{X})$ with observation densities $(g_{\theta})_{\theta\in
\sT
}$. Further, suppose that there are continuous maps $\pi\dvtx \sT\times
\tilde{\sX} \to\sX$ and $\tilde{T} \dvtx\sT\times\tilde{\sX} \to
\tilde
{\sX}$ such that:
\begin{longlist}[(iii)]
\item[(i)] for each $\theta$, we have that $\pi_{\theta} \circ
\tilde
{T}_{\theta} = T_{\theta} \circ\pi_{\theta}$;
\item[(ii)] for each $\theta$, there is a unique probability measure
$\tilde{\mu}_{\theta}$ on $\tilde{\sX}$ such that $\tilde{\mu
}_{\theta}
\circ\pi_{\theta}^{-1} = \mu_{\theta}$;
\item[(iii)] for each $\theta$, the map $\pi_{\theta}$ is injective
$\tilde{\mu}_{\theta}$-a.s.
\end{longlist}
For $x$ in $\tilde{\sX}$ and $\theta$ in $\sT$, define $\tilde
{g}_{\theta}(\cdot|x) = g_{\theta}(\cdot|\pi_{\theta}(x))$. Then
$(\tilde{T}_{\theta},\tilde{\mu}_{\theta})_{\theta\in\sT}$ is a
family of dynamical systems on $(\tilde{\sX},\tilde{\mathcal{X}})$ with
observation densities $(\tilde{g}_{\theta})_{\theta\in\sT}$. In this
situation, we say that $(T_{\theta},\mu_{\theta},g_{\theta
})_{\theta
\in\sT}$ is an isomorphic factor of $(\tilde{T}_{\theta},\tilde
{\mu
}_{\theta},\tilde{g}_{\theta})_{\theta\in\sT}$, and $\pi$ is the
factor map. The following proposition addresses the consistency of
maximum likelihood estimation for isomorphic factors. Its proof is
straightforward and omitted.

%
\begin{prop} \label{FactorProp}
Suppose that $(T_{\theta},\mu_{\theta},g_{\theta})_{\theta\in\sT}$
is an isomorphic factor of $(\tilde{T}_{\theta},\tilde{\mu}_{\theta
},\tilde{g}_{\theta})_{\theta\in\sT}$. Then maximum likelihood
estimation is consistent for $(T_{\theta},\mu_{\theta},g_{\theta
})_{\theta\in\sT}$ if and only if maximum likelihood estimation is
consistent for $(\tilde{T}_{\theta},\tilde{\mu}_{\theta},\tilde
{g}_{\theta})_{\theta\in\sT}$.
\end{prop}

For the sake of brevity, we defer precise definitions for Axiom~A
systems to Supplementary Appendix~C \cite{supplement}.

We consider families of Axiom~A systems as follows. Suppose that $f
\dvtx\sT\times\sX\to\sX$ is a parametrized family of diffeomorphisms
such that:
\begin{longlist}[(iii)]
\item[(i)] $\theta\mapsto f_{\theta}$ is H\"{o}lder continuous;
\item[(ii)] there exists $\alpha>0$ such that for each $\theta$, the
map $f_{\theta}$ is $C^{1+\alpha}$;
\item[(iii)] for each $\theta$, $\Omega(f_{\theta})$ is an Axiom~A
attractor and the restriction $f_{\theta}|_{\Omega(f_{\theta})}$ is
topologically mixing;
\item[(iv)] for each $\theta$, the measure $\mu_{\theta}$ is the unique
SRB measure corresponding to $f_{\theta}$ \cite{Bowen}, Theorem~4.1.
\end{longlist}
If these conditions are satisfied, then we say that $(f_{\theta},\mu
_{\theta})_{\theta\in\sT}$ is a parametrized family of Axiom~A
systems on $(\sX,\mathcal{X})$.

%
\begin{thmm} \label{AxiomAThm}
Suppose that $(f_{\theta},\mu_{\theta})_{\theta\in\sT}$ is a
parametrized family of Axiom~A systems on $(\sX,\mathcal{X})$. Further,
suppose that $(g_{\theta})_{\theta\in\sT}$ is a family of
observations densities satisfying the following conditions: observation
integrability \textup{(S2)} and \textup{(S3)} and
observation regularity \textup{(M3)} and \textup{(L2)}. Then maximum likelihood estimation is consistent.
\end{thmm}

The proof of Theorem~\ref{AxiomAThm} appears in the Supplementary
Appendix C \cite{supplement}.

\section{Proof of the main result} \label{SectionProofs}

Propositions \ref{ErgProp}--\ref{StrictIneqProp} are used in the proof
of Theorem~\ref{mainThm}, which is given at the end of the present section.

%
\begin{prop} \label{ErgProp}
Suppose that condition \textup{(S1)} (ergodicity) holds. Then
the process $(Y_k)$ is ergodic under $\mathbb{P}^Y_{\theta_0}$.
\end{prop}

\begin{pf}
Let $m>0$ be arbitrary, and let $A$ and $B$ be Borel subsets of $\sY
^{m+1}$. To obtain the ergodicity of $\{Y_k\}_k$, it suffices to show
that (see \cite{Petersen})
%
%
\begin{equation}
\label{eqn:ergver} \lim_n \frac{1}{n} \sum
_{k=0}^n \mathbb{P}^Y_{\theta_0}
\bigl(Y_0^m \in A, Y_k^{k+m} \in B
\bigr) = \mathbb{P}^Y_{\theta_0}\bigl(Y_0^m
\in A\bigr) \mathbb{P}^Y_{\theta
_0}\bigl(Y_0^{m}
\in B\bigr).
\end{equation}
For $x \in\sX$, define
\[
\eta_A(x) = \int\mathbf{1}_{A} \bigl(y_0^m
\bigr) p_{\theta_0} \bigl(y_0^m \mid x \bigr) \,d
\nu^{m}\bigl(y_0^m\bigr),
\]
and define $\eta_B(x)$ similarly. For $k > m$, by the conditional
independence of $Y_0^m$ and $Y_k^{k+m}$ given $\theta_0$ and $X_0=x$,
we have that
\begin{eqnarray*}
& &\mathbb{P}^Y_{\theta_0}\bigl( Y_0^m
\in A, Y_k^{k+m} \in B\bigr)
\\
&&\qquad = \int\int\mathbf{1}_{A} \bigl(y_0^m
\bigr) \mathbf{1}_B \bigl(y_k^{k+m} \bigr)
p_{\theta_0} \bigl(y_0^{n+m} \mid x \bigr) \,d\nu
^{n+m}\bigl(y_0^{n+m}\bigr) \,d\mu_{\theta_0}(x)
\\
&&\qquad = \int\biggl( \int\mathbf{1}_{A} \bigl(y_0^m
\bigr) p_{\theta
_0} \bigl(y_0^m \mid x \bigr) \,d
\nu^{m}\bigl(y_0^m\bigr)
\\
&&\hspace*{14pt}\qquad\quad{}\times\int\mathbf{1}_{B} \bigl(y_k^{k+m} \bigr)
p_{\theta_0} \bigl(y_k^{k+m} \mid T_{\theta_0}^k(x)
\bigr) \,d\nu^{m}\bigl(y_k^{k+m}\bigr) \biggr) \,d
\mu_{\theta_0}(x)
\\
&&\qquad = \int\eta_A(x) \eta_B \bigl(T_{\theta_0}^k(x)
\bigr) \,d\mu_{\theta_0}(x),
\end{eqnarray*}
where we have used Fubini's theorem.
Since $m$ is fixed, we have that
\begin{eqnarray*}
&&\lim_n \frac{1}{n}  \sum
_{k=0}^n \mathbb{P}^Y_{\theta_0}
\bigl(Y_0^m \in A, Y_k^{k+m} \in B
\bigr)
\\
&&\qquad = \lim_n \Biggl( \frac{1}{n} \sum
_{k=0}^m \mathbb{P}^Y_{\theta
_0}
\bigl(Y_0^m \in A, Y_k^{k+m} \in B
\bigr)
\\
&&\hspace*{18pt}\qquad\quad{} + \frac{1}{n} \sum_{k=m+1}^{n} \int
\eta_A(x) \eta_B \bigl(T_{\theta_0}^k(x)
\bigr) \,d\mu_{\theta_0}(x) \Biggr)
\\
&&\qquad = \lim_n \frac{1}{n} \sum
_{k=m+1}^{n} \int\eta_A(x)
\eta_B \bigl(T_{\theta_0}^k(x) \bigr) \,d
\mu_{\theta_0}(x).
\end{eqnarray*}
Since $(T_{\theta_0},\mu_{\theta_0})$ is ergodic, an alternative
characterization of ergodicity (see \cite{Petersen}) gives that
\begin{eqnarray*}
\lim_n \frac{1}{n} \sum_{k=0}^n
\mathbb{P}^Y_{\theta_0} \bigl(Y_0^m \in
A, Y_k^{k+m} \in B \bigr) & = &\lim_n
\frac{1}{n} \sum_{k=m+1}^{n} \int\eta
_A(x) \eta_B \bigl(T_{\theta_0}^k(x)
\bigr) \,d\mu_{\theta_0}(x)
\\
& = &\int\eta_A(x) \,d\mu_{\theta_0}(x) \int\eta_B(x)
\,d\mu_{\theta
_0}(x)
\\
& = &\mathbb{P}^Y_{\theta_0}\bigl(Y_0^m
\in A\bigr) \mathbb{P}^Y_{\theta
_0}\bigl(Y_0^m
\in B\bigr).
\end{eqnarray*}
Thus we have verified equation \eqref{eqn:ergver}, and the proof is complete.
\end{pf}

For the following propositions, recall our notation that
\[
\gamma_{\theta}(y) = \sup_x g_{\theta}(y \mid
x).
\]

%
\begin{prop} \label{gSMB}
Suppose that conditions \textup{(S1)} and \textup{(S2)} hold. Then there exists $h(\theta_0) \in(-\infty
,\infty)$ such that
\[
h(\theta_0) = \lim_n \E_{\theta_0} \biggl(
\frac{1}{n} \log p_{\theta
_0}\bigl(Y_0^n\bigr)
\biggr).
\]
Moreover, the following equality holds $\P_{\theta_0}$-a.s.:
\[
h(\theta_0) = \lim_n \frac{1}{n} \log
p_{\theta_0}\bigl(Y_0^n\bigr).
\]
\end{prop}

\begin{pf}
The proposition is a direct application of Barron's generalized
Shannon--McMillan--Breiman theorem \cite{Barron1985}. Here we simply
check that the hypotheses of that theorem hold in our setting. Since
condition (S1) (ergodicity) holds, Proposition~\ref
{ErgProp} gives $(Y_k)$ is stationary and ergodic under $\P_{\theta
_0}$. By definition, $Y_0^n$ has density $p_{\theta_0}(Y_0^n)$ with
respect to the $\sigma$-finite measure $\nu^{n}$. The measure $\nu^n$
is a product of the measure $\nu$ taken $n+1$ times. As such, the
sequence $\{\nu^n\}$ clearly satisfies Barron's condition that this
sequence is ``Markov with stationary transitions.'' Define $D_n = \E
_{\theta_0}( \log p_{\theta_0}(Y_0^{n+1})) - \E_{\theta_0}( \log
p_{\theta_0}(Y_0^{n}))$. Let us show that for $n >0$, we have that
%
%
\begin{equation}
\label{FiniteEntropy} \E_{\theta_0}\bigl( \bigl|\log p_{\theta_0}
\bigl(Y_0^n\bigr)\bigr| \bigr) < \infty,
\end{equation}
which clearly implies that $- \infty< D_n < \infty$. Once \eqref
{FiniteEntropy} is established, we will have verified all of the
hypotheses of Barron's generalized Shannon--McMillan--Breiman theorem,
and the proof of the proposition will be complete.

Observe that the first part of the integrability condition (S2) gives that
%
%
\begin{equation}
\label{Eqn:PosPartFinite} \E_{\theta_0} \bigl[ \log^+ p_{\theta_0}
\bigl(Y_0^n\bigr) \bigr] \leq(n+1) \E_{\theta_0}
\bigl[ \log^+ \gamma_{\theta_0}(Y_0) \bigr] < \infty.
\end{equation}
Then the second part of the integrability condition (S2) implies that
%
%
\begin{eqnarray}
\label{Eqn:NegPartFinite} %
\E_{\theta_0} \bigl[ \log p_{\theta_0}
\bigl(Y_0^n\bigr) \bigr] & = &\E_{\theta
_0} \biggl[
\log\frac{ p_{\theta_0}(Y_0^n) }{\prod_{k=0}^n \int
g_{\theta_0}(Y_k \mid x) \,d\mu_{\theta_0}(x)} \biggr]
\nonumber\\
&&{} + \E_{\theta_0} \Biggl[ \sum_{k=0}^n
\log\int g_{\theta
_0}(Y_k \mid x) \,d\mu_{\theta_0}(x) \Biggr]
\nonumber
\\[-8pt]
\\[-8pt]
\nonumber
& \geq&-(n+1) \E_{\theta_0} \biggl[ \biggl| \log\int g_{\theta_0}(Y_0
\mid x) \,d\mu_{\theta_0}(x) \biggr| \biggr]
\\
& > &-\infty,
\nonumber
\end{eqnarray}
where we have used that relative entropy is nonnegative. By \eqref
{Eqn:PosPartFinite} and \eqref{Eqn:NegPartFinite}, we conclude that
\eqref{FiniteEntropy} holds, which completes the proof.
\end{pf}

The following proposition is used in the proof of Theorem~\ref{mainThm}
to given an almost sure bound for the normalized log-likelihoods in
terms of quantities involving only expectations.

%
\begin{prop} \label{BasicEstimate}
Suppose that conditions \textup{(S1)}, \textup{(S3)}
and \textup{(S5)} hold. Let $\ell$ be as in condition \textup{(S5)}.
Then for $\theta' \notin[\theta_0]$, there exists a
neighborhood $U$ of $\theta'$ such that for each $m >0$, the following
inequality holds $\P_{\theta_0}$-a.s.:
\begin{eqnarray*}
\limsup_{n \to\infty} \sup_{\theta\in U } \frac{1}{n}
\log p_{\theta}\bigl(Y_0^n\bigr) &\leq&
\frac{1}{m+\ell} \E_{\theta_0} \Bigl( \sup_{\theta\in U} \log
p_{\theta}\bigl(Y_0^m\bigr) \Bigr)
\\
&&{} + \frac{\ell}{m+\ell} \E_{\theta_0} \Bigl( \sup_{\theta\in
U}
\log^+ \gamma_{\theta}(Y_0) \Bigr)
\\
&&{} + \frac{1}{m+\ell} \E_{\theta_0} \Bigl( \sup_{\theta\in U} \log
C_m\bigl(\theta,Y_0^m\bigr) \Bigr).
\end{eqnarray*}
\end{prop}

Informally, in the proof of Proposition~\ref{BasicEstimate}, we use the
mixing property from condition (S5) to parse a
sequence of observations into alternating sequences of ``large blocks''
and ``small blocks,'' and then the ergodicity and integrability
conditions finish the proof. More specifically, we break up the
sequence of observations $Y_0^n$ into alternating blocks of length $m$
and $\ell$, where $\ell$ is given by condition~(S5).
\begin{pf}
Let $\theta' \notin[\theta_0]$. Fix a neighborhood $U$ of $\theta'$
so that the conclusions of both condition (S3) and
condition (S5) hold. Let $m>0$ be arbitrary, and let
$\ell$ be as in condition (S5). We consider sequences
of observations of length $n$, where $n$ is a large integer. These
sequences of observations will be parsed into alternating blocks of
lengths $m$ and $\ell$, respectively, starting from an offset of size
$s$ and possibly ending with a remainder sequence. For the sake of
notation, we use interval notation to denote intervals of integers. For
$n > 2(m+\ell)$ and $s$ in $[0,m+\ell)$, let $R = R(s,m,\ell,n) \in
[0,m+\ell)$ and $k = k(s,m,\ell,n) \geq0$ be defined by the condition
$n = s+ k(m+\ell)+R$. Then we partition $[0,n]$ as follows:
\begin{eqnarray*}
B_s & =& [0, s ),
\\
I_s(j) & = &\bigl[s+(m+\ell) (j-1), s+(m+\ell) (j-1) + m \bigr)\qquad \mbox{for
} 1 \leq j \leq k,
\\
J_s(j) & = &\bigl[s+(m+\ell) (j-1)+m, s+(m+\ell)j \bigr)\qquad\mbox{for } 1
\leq j \leq k,
\\
E_s & =& \bigl[s+t(m+\ell), n \bigr].
\end{eqnarray*}
Given a sequence $Y_0^n$ of observations, we define the following
subsequences of $Y_0^n$ according to the above partitions of $[0,n]$:
\begin{eqnarray*}
b_s &=& Y|_{B_s},
\\
w_s(j) & =& Y|_{I_s(j)}\qquad \mbox{for } 1 \leq j \leq k,
\\
v_s(j) & =& Y|_{J_s(j)}\qquad \mbox{for } 1 \leq j \leq k,
\\
e_s & =& Y|_{E_s}.
\end{eqnarray*}
For a sequence $y_0^t$ in $\sY^{t+1}$, define
\[
\gamma_{\theta}\bigl(y_0^t\bigr) = \prod
_{j=0}^t \gamma_{\theta}(y_j) =
\prod_{j=0}^t \sup_x
g_{\theta}(y_j \mid x).
\]
Then for $\theta$ in $U$, it follows from condition (S5) that
\[
p_{\theta}\bigl(Y_0^{n}\bigr) \leq
\gamma_{\theta}(b_s) \gamma_{\theta}(e_s)
\cdot\prod_{j=1}^k \gamma_{\theta}
\bigl(v_s(j)\bigr) \cdot\prod_{j=1}^{k}C_{m}
\bigl(\theta,w_s(j)\bigr) \cdot\prod_{j=1}^{k}
p_{\theta}\bigl(w_s(j)\bigr).
\]
Taking the logarithm of both sides and averaging over $s$ in $[0,
m+\ell
)$, we obtain
%
%
\begin{eqnarray}
\label{Eqn:ThreePieces} %
\log p_{\theta}\bigl(Y_0^{n}
\bigr) &\leq& \frac{1}{m+\ell} \sum_{s=0}^{m+\ell-1}
\sum_{j=1}^{k} \bigl[ \log p_{\theta}
\bigl(w_s(j)\bigr) + \log C_{m}\bigl(\theta
,w_s(j)\bigr) \bigr]\nonumber
\\
&&{} + \frac{1}{m+\ell} \sum_{s=0}^{m+\ell-1} \sum
_{j=1}^k \log\gamma_{\theta}
\bigl(v_s(j)\bigr)
\\
&&{} + \frac{1}{m+\ell} \sum_{s=0}^{m+\ell-1}
\bigl[ \log\gamma_{\theta
}(b_s) + \log\gamma_{\theta}(e_s)
\bigr].
\nonumber
\end{eqnarray}
Let us now take the supremum over $\theta$ in $U$ in \eqref
{Eqn:ThreePieces} and evaluate the limits of the three terms on the
right-hand side as $n$ tends to infinity.

Let $\xi_1 \dvtx\sY^{m+1} \to\R$ and $\xi_2 \dvtx\sY^{m+1} \to\R$ be
defined by
\begin{eqnarray*}
\xi_1\bigl(y_0^m\bigr)& =& \sup
_{\theta\in U} \log p_{\theta}\bigl(y_0^m
\bigr),\\
 \xi_2\bigl(y_0^m\bigr) &=& \sup
_{\theta\in U} \log C_m\bigl(\theta,y_0^m
\bigr).
\end{eqnarray*}
With this notation, we have that
\begin{eqnarray*}
&&\frac{1}{n} \sum_{s=0}^{m+\ell-1}  \sum
_{j=1}^{k} \Bigl[ \sup_{\theta\in U}
\log p_{\theta}\bigl(w_s(j)\bigr) + \sup_{\theta\in U}
\log C_{m}\bigl(\theta,w_s(j)\bigr) \Bigr]
\\
&&\qquad = \frac{1}{n} \sum_{i=0}^n \bigl[
\xi_1\bigl(Y_i^{i+m}\bigr) + \xi
_2\bigl(Y_i^{i+m}\bigr) \bigr].
\end{eqnarray*}
Since $(Y_k)$ is ergodic (by Proposition~\ref{ErgProp}), it follows
from Birkhoff's ergodic theorem and conditions (S3)
and (S5) that the following limit exists $\P
_{\theta
_0}$-a.s.:
%
%
\begin{eqnarray}
\label{Term1} %
&&\lim_n \frac{1}{n} \sum
_{s=0}^{m+\ell-1}  \sum
_{j=1}^{k} \Bigl[ \sup_{\theta\in U} \log
p_{\theta}\bigl(w_s(j)\bigr) + \sup_{\theta\in U}
\log C_{m}\bigl(\theta,w_s(j)\bigr) \Bigr]
\nonumber\\
&&\qquad = \lim_n \frac{1}{n} \sum
_{i=0}^n \bigl[ \xi_1
\bigl(Y_i^{i+m}\bigr) + \xi_2
\bigl(Y_i^{i+m}\bigr) \bigr]
\nonumber
\\
& &\qquad= \E_{\theta_0} \bigl[ \xi_1\bigl(Y_0^m
\bigr) \bigr] + \E_{\theta_0} \bigl[ \xi_2\bigl(Y_0^m
\bigr) \bigr]
\\
&&\qquad = \E_{\theta_0} \Bigl[ \sup_{\theta\in U} \log p_{\theta}
\bigl(Y_0^m\bigr) \Bigr]\nonumber\\
&&\qquad\quad{} + \E_{\theta_0} \Bigl[ \sup
_{\theta\in U} \log C_{m}\bigl(\theta,Y_0^m
\bigr) \Bigr].
\nonumber
\end{eqnarray}

Similarly, using Birkhoff's ergodic theorem and condition (S3), we have that the following holds $\P_{\theta_0}$-a.s.:
\begin{eqnarray}
\label{Term2} %
\limsup_n \frac{1}{n} \sum
_{s=0}^{m+\ell-1} \sum
_{j=1}^{k} \sup_{\theta\in U} \log
\gamma_{\theta}\bigl(v_s(j)\bigr) & \leq&\limsup
_n \frac
{1}{n} \sum_{i=0}^n
\sup_{\theta\in U} \log^+\gamma_{\theta
}\bigl(Y_i^{i+\ell-1}
\bigr)
\nonumber
\\
& \leq&\ell\limsup_n \frac{1}{n} \sum
_{i=0}^n \sup_{\theta\in U} \log^+
\gamma_{\theta}(Y_i)
\\
& =& \ell\E_{\theta_0} \Bigl[ \sup_{\theta\in U} \log^+ \gamma
_{\theta}(Y_0) \Bigr].
\nonumber
\end{eqnarray}
Finally, Birkhoff's ergodic theorem and condition (S3) again imply that the following limit holds $\P
_{\theta
_0}$-a.s.:
%
%
\begin{equation}
\label{Term3} \lim_n \frac{1}{n} \sum
_{s=0}^{m+\ell-1} \Bigl[ \sup_{\theta\in U} \log^+
\gamma_{\theta}(b_s) + \sup_{\theta\in U} \log^+
\gamma_{\theta
}(e_s) \Bigr] = 0,
\end{equation}
where we have used that $\max(|B_s|,|E_s|) \leq m+\ell$.

Combining the inequalities in (\ref{Eqn:ThreePieces})--(\ref{Term3}),
we obtain that
\begin{eqnarray*}
\limsup_{n \to\infty} \sup_{\theta\in U} \frac{1}{n}
\log p_{\theta
}\bigl(Y_0^{n}\bigr) &\leq&
\frac{1}{m+\ell} \E_{\theta_0} \Bigl[ \sup_{\theta
\in U} \log
p_{\theta}\bigl(Y_0^m\bigr) \Bigr]
\\
&&{} + \frac{1}{m+\ell} \E_{\theta_0} \Bigl[ \sup_{\theta\in U} \log
C_{m}\bigl(\theta,Y_0^m\bigr) \Bigr]
\\
&&{} + \frac{\ell}{m+\ell} \E_{\theta_0} \Bigl[ \sup_{\theta\in
U}
\log^+ \gamma_{\theta}(Y_0) \Bigr],
\end{eqnarray*}
as desired.
\end{pf}

The following proposition is a direct application of Lemma~10 in \cite
{DMOvH2011} 
to the present setting, and we omit the proof.
%

\begin{prop} \label{IdentifiabilityProp}
Suppose that the following conditions hold: ergodicity~\textup{(S1)},
logarithmic integrability at $\theta_0$ \textup{(S2)} and exponential
identifiability \textup{(S6)}. Then for $\theta\notin[\theta_0]$, it holds that
\[
\limsup_n \frac{1}{n} \E_{\theta_0} \bigl[ \log
p_{\theta}\bigl(Y_0^{n}\bigr) \bigr] < h(
\theta_0).
\]
\end{prop}

The following proposition provides an essential estimate in the proof
of Theorem~\ref{mainThm}.

%
\begin{prop} \label{StrictIneqProp}
Suppose that conditions \textup{(S1)--(S6)}
hold, and let $\ell$ be as in~\textup{(S5)}. Then for
$\theta'
\notin[\theta_0]$, there exists $m>0$ and a neighborhood $U$ of
$\theta
'$ such that
\begin{eqnarray*}
h(\theta_0) &> & \frac{1}{m+\ell} \E_{\theta_0} \Bigl[ \sup
_{\theta\in
U} \log p_{\theta}\bigl(Y_0^m
\bigr) \Bigr]
\\
&&{} + \frac{\ell}{m+\ell} \E_{\theta_0} \Bigl[ \sup_{\theta\in
U}
\log^+ \gamma_{\theta}(Y_0) \Bigr]
\\
&&{} + \frac{1}{m+\ell} \E_{\theta_0} \Bigl[ \sup_{\theta\in U} \log
C_{m}\bigl(\theta',Y_0^{m}\bigr)
\Bigr].
\end{eqnarray*}
\end{prop}

\begin{pf}
Suppose $\theta' \notin[\theta_0]$. By Proposition~\ref
{IdentifiabilityProp}, there exists $\varepsilon>0$ such that
%
%
\begin{equation}
\label{Eqn:Estimate1} \limsup_n \frac{1}{n}
\E_{\theta_0} \bigl[ \log p_{\theta'}\bigl(Y_0^{n}
\bigr) \bigr] < h(\theta_0)-\varepsilon.
\end{equation}
By conditions (S3) (logarithmic integrability away
from $\theta_0$) and (S5) (mixing), there exists a
neighborhood $U'$ of $\theta'$ and $m_0>0$ such that for $m \geq m_0$,
we have that
%
%
\begin{eqnarray}
\label{Eqn:Estimate2} %
\varepsilon/2& > & \frac{\ell}{m+\ell} \E_{\theta_0}
\Bigl[ \sup_{\theta
\in U'} \log^+ \gamma_{\theta}(Y_0)
\Bigr]
\nonumber
\\[-8pt]
\\[-8pt]
\nonumber
&&{} + \frac{1}{m+\ell} \E_{\theta_0} \Bigl[ \sup_{\theta\in U'} \log
C_{m}\bigl(\theta,Y_0^{m}\bigr) \Bigr].
\nonumber
\end{eqnarray}
Fix $m \geq m_0$ such that
%
%
\begin{equation}
\label{Eqn:Estimate3} \frac{1}{m+\ell} \E_{\theta_0} \bigl[ \log
p_{\theta'}
\bigl(Y_0^{m}\bigr) \bigr] < \limsup_n
\frac{1}{n} \E_{\theta_0} \bigl[ \log p_{\theta
'}
\bigl(Y_0^{n}\bigr) \bigr] + \varepsilon/4.
\end{equation}
For $\eta>0$, let $B(\theta',\eta)$ denote the ball of radius $\eta$
about $\theta'$ in $\sT$. For $\eta$ such that $B(\theta',\eta)
\subset
U'$, we have that
\[
\sup_{\theta\in B(\theta',\eta)} \log p_{\theta}\bigl(Y_0^m
\bigr) \leq\sum_{k=0}^m \sup
_{\theta\in U'} \log^+ \gamma_{\theta}(Y_k).
\]
The sum above is integrable with respect to $\P_{\theta_0}$ and does
not depend on $\eta$. Then (the reverse) Fatou's Lemma implies that
\[
\limsup_{\eta\to0} \E_{\theta_0} \Bigl[ \sup
_{\theta\in
B(\theta
',\eta)} \log p_{\theta}\bigl(Y_0^{m}
\bigr) \Bigr] \leq\E_{\theta_0} \Bigl[ \limsup_{\eta\to0} \sup
_{\theta\in B(\theta',\eta)} \log p_{\theta
}\bigl(Y_0^{m}
\bigr) \Bigr].
\]
By condition (S4) [upper semi-continuity of $\theta\mapsto
p_{\theta}(Y_0^m)$], we see that
\[
\E_{\theta_0} \Bigl[ \limsup_{\eta\to0} \sup
_{\theta\in
B(\theta
',\eta)} \log p_{\theta}\bigl(Y_0^{m}
\bigr) \Bigr] \leq\E_{\theta_0} \bigl[ \log p_{\theta'}
\bigl(Y_0^{m}\bigr) \bigr].
\]
Now by an appropriate choice of $\eta>0$, we have shown that there
exists a neighborhood $U \subset U'$ of $\theta'$ such that
%
%
\begin{equation}
\label{Eqn:Estimate4} \frac{1}{m+\ell} \E_{\theta_0} \Bigl[ \sup
_{\theta\in U} \log p_{\theta}\bigl(Y_0^{m}
\bigr) \Bigr] < \frac{1}{m+\ell} \E_{\theta_0} \bigl[ \log p_{\theta'}
\bigl(Y_0^{m}\bigr) \bigr] + \varepsilon/4.
\end{equation}
Combining estimates \eqref{Eqn:Estimate1}--\eqref{Eqn:Estimate4}, we
obtain the desired inequality.
\end{pf}

\begin{pf*}{Proof of Theorem~\ref{mainThm}}
Let $h(\theta_0)$ be defined as in Proposition~\ref{gSMB}. We prove the
theorem by showing the following statement: for each closed set $C$ in
$\sT$ such that $C \cap[\theta_0] = \varnothing$, it holds that
%
%
\begin{equation}
\label{suffIneq1} \limsup_n \sup_{\theta\in C}
\frac{1}{n} \log p_{\theta}\bigl(Y_0^{n}\bigr)
< h(\theta_0).
\end{equation}

Let $C$ be a closed subset of $\sT$ such that $C \cap[\theta_0] =
\varnothing$. Since $\sT$ is compact, $C$ is compact. Suppose that for
each $\theta' \in C$, there exists a neighborhood $U$ of $\theta'$
such that
%
%
\begin{equation}
\label{suffIneq2} \limsup_n \sup_{\theta\in U \cap C}
\frac{1}{n} \log p_{\theta
}\bigl(Y_0^{n}\bigr)
< h(\theta_0).
\end{equation}
Then by compactness, we would conclude that (\ref{suffIneq1}) holds and
thus complete the proof of the theorem.

Let $\theta'$ be in $C$. Let us now show that there exists a
neighborhood $U$ of $\theta'$ such that (\ref{suffIneq2}) holds. Since
$\theta'$ is in $C$, we have that $\theta' \notin[\theta_0]$. Let
$\ell
$ be as in (S5). By Proposition~\ref{StrictIneqProp},
there exists $m >0$ and a neighborhood $U'$ of $\theta'$ such that
%
%
\begin{eqnarray}
\label{Eqn:MainProofEst1} %
h(\theta_0) &> & \frac{1}{m+\ell}
\E_{\theta_0} \Bigl[ \sup_{\theta\in
U'} \log p_{\theta}
\bigl(Y_0^m\bigr) \Bigr]
\nonumber\\
& &{}+ \frac{\ell}{m+\ell} \E_{\theta_0} \Bigl[ \sup_{\theta\in
U'}
\log^+ \gamma_{\theta}(Y_0) \Bigr]
\\
& &{}+ \frac{1}{m+\ell} \E_{\theta_0} \Bigl[ \sup_{\theta\in U'} \log
C_{m}\bigl(\theta,Y_0^{m}\bigr) \Bigr].
\nonumber
\end{eqnarray}
By Proposition~\ref{BasicEstimate}, there exists a neighborhood $U
\subset U'$ of $\theta'$ such that
%
%
\begin{eqnarray}
\label{Eqn:MainProofEst2} %
\limsup_{n \to\infty} \sup
_{\theta\in U } \frac{1}{n} \log p_{\theta}
\bigl(Y_0^n\bigr) &\leq& \frac{1}{m+\ell}
\E_{\theta_0} \Bigl[ \sup_{\theta\in U} \log p_{\theta}
\bigl(Y_0^m\bigr) \Bigr]\nonumber
\\
&&{} + \frac{\ell}{m+\ell} \E_{\theta_0} \Bigl[ \sup_{\theta\in
U}
\log^+ \gamma_{\theta}(Y_0) \Bigr]
\\
&&{} + \frac{1}{m+\ell} \E_{\theta_0} \Bigl[ \sup_{\theta\in U} \log
C_{m}\bigl(\theta',Y_0^{m}\bigr)
\Bigr].
\nonumber
\end{eqnarray}
Combining \eqref{Eqn:MainProofEst1} and \eqref{Eqn:MainProofEst2}, we
obtain (\ref{suffIneq2}), which completes the proof of the theorem.
\end{pf*}


\section{Concluding remarks} \label{SectionConclusion}

In this paper, we demonstrate how the properties of a family of
dynamical systems affect the asymptotic consistency of maximum
likelihood parameter estimation. We have exhibited a collection of
general statistical conditions on families of dynamical systems
observed with noise, and we have shown that under these general
conditions, maximum likelihood estimation is a consistent method of
parameter estimation. Furthermore, we have shown that these general
conditions are indeed satisfied by some classes of well-studied
families of dynamical systems. As mentioned in the \hyperref
[Introduction]{Introduction}, our
results can be considered as a theoretical validation of the notion
from dynamical systems that these classes of systems have ``good''
statistical properties.

However, there remain interesting families of systems to which our
results do not apply, including some classes of systems that are also
believed to have ``good'' statistical properties. In particular, the
class of systems modeled by Young towers with exponential tail \cite
{Young1998} has exponential decay of correlations and certain large
deviations estimates \cite{RBY2008}. These families include a positive
measure set of maps from the quadratic family [$x \mapsto ax(1-x)$] and
the H\'{e}non family, as well as certain billiards and many other
systems of physical and mathematical interest \cite{Young1998}. In
short, the setting of systems modeled by Young towers with exponential
tail provides a very attractive setting in which to consider
consistency of maximum likelihood estimation. Unfortunately, our proof
does not apply to systems in this setting in general, mainly due to the
presence of the mixing condition (S5), which is not
satisfied by these systems in general.

A natural next step might be to obtain rates of convergence and derive
central limit theorems for maximum likelihood estimation. To this end,
it might be possible to build off of analogous results for HMMs \cite
{BRR1998,JensenPetersen1999}. We leave these questions for future work.


\begin{supplement}[id=suppA]
\stitle{Supplement to ``Consistency of maximum likelihood estimation
for some dynamical systems''}
\slink[doi]{10.1214/14-AOS1259SUPP} 
\sdatatype{.pdf}
\sfilename{aos1259\_supp.pdf}
\sdescription{We provide three technical appendices. In Appendix~A, we present proofs of Propositions \ref
{USClemma2}, \ref{SubaddPropMain} and \ref
{IdentifiabilityFromLargeDeviations}. In Appendix~B, we discuss shifts of finite type and Gibbs measures
and prove Theorem~\ref{GibbsFamiliesThm}. Finally, Appendix~C contains definitions for Axiom~A systems, as well as
a proof of Theorem~\ref{AxiomAThm}.}
\end{supplement}

\printaddresses
\end{document}